\input amstex

\define\R{\text{\rm{Re}}\,}
\define\I{\text{\rm{Im}}\,}

\documentstyle {amsppt}
\magnification=\magstep1 \NoRunningHeads

\topmatter
\title
 The convergence of a sequence of polynomials and the distribution of
their zeros
 \endtitle
\author
Min-Hee Kim, Young-One Kim and Jungseob Lee
\endauthor


\abstract Suppose that $\langle f_n \rangle$ is a sequence of
polynomials, $\langle f_n^{(k)}(0)\rangle$ converges for every
non-negative integer $k$, and that the limit is not $0$ for some
$k$. It is shown that if all the zeros of $f_1, f_2, \dots$ lie in
the closed upper half plane $\I z\geq 0$, or if $f_1, f_2, \dots$
are real polynomials and the numbers of their non-real zeros are
uniformly bounded, then the sequence converges uniformly on compact
sets in the complex plane. The results imply a theorem of Benz and a
conjecture of P\'olya.
\endabstract

\dedicatory Dedicated to the memory of the late Professor Jehpill
Kim (1930--2016).
\enddedicatory

\address
Department of Mathematical Sciences, Seoul National University,
Seoul  08826,    Korea
\endaddress

\email alsgml01\@snu.ac.kr
\endemail

\address
Department of Mathematical Sciences and Research Institute of
Mathematics, Seoul National University, Seoul  08826,  Korea
\endaddress

\email kimyone\@snu.ac.kr
\endemail

\address
Department of Mathematics, Ajou University, Suwon 16499,  Korea
\endaddress

\email jslee\@ajou.ac.kr
\endemail



\thanks
This work was supported by the SNU Mathematical Sciences Division
for Creative Human Resources Development.
\endthanks

\endtopmatter

\document

\heading 1. Introduction \endheading

This paper is concerned with the distribution of zeros of entire
functions and the convergence of polynomial sequences. Let $f$ be an
entire function which is not identically equal to $0$. We denote the
zero set of $f$ by $\Cal Z(f)$, that is, $ \Cal Z(f)=\{z\in\Bbb C :
f(z)=0\}$. If $X$ is a subset of the complex plane, the number of
zeros of $f$ that lie in $X$ counted according to their
multiplicities will be denoted by $N(f;X)$.  In the case where $f$
is identically equal to $0$, we set $\Cal Z(f)=\emptyset$ and
$N(f;X)=0$. The closed upper half plane $\{z\in\Bbb C : \I z\geq
0\}$ is denoted by $\Bbb H$, and for $r>0$ the open disk $\{z\in\Bbb
C : |z|<r\}$ is denoted by $D_r$.

The {\it P\'olya-Obrechkoff class} $\Cal{PO}$ is the collection of
entire functions $f$ which can be represented in the form
$$
f(z)=cz^m e^{-\alpha z^2 + \beta
z}\prod_j\left(1-\frac{z}{a_j}\right)e^{z/a_j} \qquad (z\in\Bbb C),
\tag 1.1
$$
where $c$ is a constant, $m$ is a nonnegative integer, $\alpha\geq
0$, $ a_1, a_2, \dots \in{\Bbb H}\setminus\{0\}$,
$\sum_{j}|a_j|^{-2}<\infty$ and $\sum_{j }\I
\left(-a_j^{-1}\right)\leq \I \beta$. If $f\in\Cal{PO}$, then it is
easy to see that there is a sequence $\langle f_n\rangle$ of
polynomials such that $f_n\to f$ uniformly on compact sets in the
complex plane and $\Cal Z (f_n)\subset {\Bbb H}$ for all $n$. In
1914 Lindwart and P\'olya proved a strong version of the converse
\cite{4, Satz II B}.

\proclaim{Theorem (Lindwart-P\'olya)} Suppose that  $\langle
f_n\rangle$ is a sequence of polynomials,  $\langle f_n\rangle$
converges uniformly on a disk $D_r$, the limit is not identically
equal to $0$, and that $\Cal Z (f_n)\subset {\Bbb H}$ for all $n$.
Then it converges to an entire function $f\in\Cal{PO}$ uniformly on
compact sets in the complex plane.
\endproclaim

Since then ``a disk $D_r$'' in the theorem has been replaced with
considerably smaller sets by several authors: Korevaar and  Loewner
replaced it with a certain arc \cite{3},  Levin replaced it with a
certain countable set, and later  Clunie and  Kuijlaars proved
Levin's theorem in a new and simple way, and extended it \cite{2}.

In this paper, we approach the problem in a different direction. If
a sequence $\langle f_n\rangle$ of analytic functions  converges
uniformly on a disk  $D_r$, then for every nonnegative integer $k$
the sequence $\langle f_n^{(k)}(0)\rangle$ of complex numbers
converges, but the converse does not hold in general. Let us say
that $\langle f_n\rangle$ converges {\it weakly} if $\langle
f_n^{(k)}(0)\rangle$ converges for every nonnegative integer $k$. In
this case, the limit is the formal power series whose $k$-th
coefficient is given by
$$
\frac{1}{k!}\lim_{n\to\infty}f_n^{(k)}(0)
$$
for every $k$.  For notational simplicity, we say that a sequence of
entire functions converges {\it strongly} if it converges uniformly
on compact sets in the complex plane.

In 1934, Benz proved that a formal power series $f$ represents an
entire function (the radius of convergence is $\infty$) in
$\Cal{PO}$ if there is a sequence $\langle f_n\rangle$ of
polynomials such that $f_n\to f$ weakly and $\Cal Z(f_n)\subset
{\Bbb H}$ for all $n$ \cite{1, Satz 2}. However, the result  implies
nothing about uniform convergence of the sequence, except in the
case of the Jensen sequence introduced in Section 3 below.

As our first result, we improve the Lindwart-P\'olya theorem as well
as Benz's result.

\proclaim{Theorem 1.1} Suppose that  $\langle f_n\rangle$ is a
sequence of polynomials, $f_n\to f$ weakly, some coefficient of $f$
is not $0$, and that $\Cal Z(f_n)\subset {\Bbb H}$ for all $n$. Then
$f$ represents an entire function in the P\'olya-Obrechkoff class
and $f_n\to f$ strongly.
\endproclaim

\remark{Remarks} (i)  As the example $f_n(z)=z^n$ shows, the
condition that {\it some coefficient of $f$ is not $0$} is
necessary. (ii) It may be remarked that in the general case any of
the convergence conditions considered in \cite{2} and  \cite{3} does
not imply weak convergence and vice versa.
\endremark

 If $f$ is of the form (1.1) with $c, \beta,a_1, a_2,
\dots\in\Bbb R$, then $f$ is said to be in the {\it Laguerre-P\'olya
class} $\Cal{LP}$. It is clear that $\Cal{LP}\subset\Cal{PO}$, $\Cal
Z(f)\subset \Bbb R$ for all $f\in\Cal{LP}$, and that $f\in\Cal{LP}$
if and only if $f\in\Cal{PO}$ and $f$ is a {\it real entire
function}, that is, $f(\Bbb R)\subset \Bbb R$. In  \cite{5}, P\'olya
introduced another extension of $\Cal{LP}$, namely the class
$\Cal{LP}^*$ of real entire functions $f$ which are of the form
$f=Pg$ where $P$ is a real polynomial and $g\in\Cal{LP}$.  If
$f\in\Cal{LP}^*$, then $N(f; \Bbb C\setminus\Bbb R)<\infty$ and
there is a sequence $\langle f_n\rangle$ of real polynomials such
that $f_n\to f$ strongly and $ N(f_n; \Bbb C\setminus\Bbb R)= N(f;
\Bbb C\setminus\Bbb R)$ for all $n$. As P\'olya stated in the same
paper,  if a sequence  $\langle f_n\rangle$ of real polynomials
converges uniformly on a disk $D_r$, the limit is not identically
equal to $0$  and $\langle N(f_n; \Bbb C\setminus\Bbb R)\rangle$ is
bounded above, then it converges strongly to an entire function in
$\Cal{LP}^*$ \cite{5, \S 4 Satz II}.  It is not clear whether the
results of \cite{2} and \cite{3} extend to the case of $\Cal{LP}^*$,
but weak convergence implies strong convergence in this case too.

\proclaim{Theorem 1.2} Suppose that  $\langle f_n\rangle$ is a
sequence of real polynomials, $f_n\to f$ weakly, some coefficient of
$f$ is not $0$, and that $\langle N(f_n; \Bbb C\setminus\Bbb
R)\rangle$ is bounded. Then $f$ represents an entire function in
$\Cal{LP}^*$  and $f_n\to f$ strongly.
\endproclaim

As we shall see in the sequel, this theorem implies a conjecture of
P\'olya which has  remained open since 1915.

In Section 2, we  prove Theorems 1.1 and 1.2 by generalizing a
theorem of Lindwart and P\'olya. Finally, we apply the results to
give a simple proof of the original  version of Benz's theorem
mentioned above and to improve some classical theorems of P\'olya
(Section 3).

\heading 2. Proofs of Theorems 1.1 and 1.2
\endheading

In this section, we obtain some generalizations (Theorems 2.3 and
2.4 below) of results in \cite{4}, and prove Theorems 1.1 and 1.2.

Let $f$ be an entire function such that $f(0)\ne 0$, and let $a_1,
a_2, \dots$ be the zeros of $f$ listed according to their
multiplicities. If $X\subset\Bbb C$ and $k$ is a positive integer,
we put
$$
\tilde s_k(f; X)=\sum_{a_j\in X} |a_j|^{-k};
$$
and in the case where $\tilde s_k(f; X)<\infty$, we put
$$
s_k(f; X)=\sum_{a_j\in X} a_j^{-k} .
$$
We also write $s_k(f)=s_k(f; \Bbb C)$ and  $\tilde s_k(f)=\tilde
s_k(f; \Bbb C)$.

Suppose that  $f$ is a polynomial and $f(0)\ne 0$. Then there is a
positive constant $r$ such that
$$
f(z)=f(0)\exp\left(-\sum_{k=1}^{\infty}\frac{1}{k}s_k(f)z^k\right)\qquad
(z\in D_r).
$$
Thus we have
$$
s_k(f)=\frac{-1}{(k-1)!}\left(\frac{f'}{f}\right)^{(k-1)}
\!\!\!\!(0)\qquad (k=1, 2, \dots),
$$
and it follows that $s_k(f)$ is a rational function of $f(0), f'(0),
\dots, f^{(k)}(0)$ for $k=1, 2, \dots$. For instance,
$s_1(f)=-f'(0)/f(0)$, $s_2(f)=\left(f'(0)^2 -
f(0)f''(0)\right)/f(0)^2$, and so on. This observation leads to the
following:

\proclaim{Proposition 2.1} Suppose that $\langle f_n \rangle$ is a
  sequence of polynomials, $f_n(0)\ne 0$ for all $n$, $\lim_{n\to\infty}f_n(0)\ne 0$ and $f_n\to f$ weakly.
  Then the sequences $\langle s_k(f_n)\rangle$, $k=1, 2, \dots$,
are all convergent. If $f$ represents an analytic function in a
neighborhood of $0$, then we have
$$
\lim_{n\to\infty}
s_k(f_n)=\frac{-1}{(k-1)!}\left(\frac{f'}{f}\right)^{(k-1)}\!\!\!\!(0)\qquad
(k=1, 2, \dots).
$$
\endproclaim

It is obvious that $|1-z|\leq e^{|z|}$ for all $z\in \Bbb C$. More
generally, we have the following, which is also trivially proved.

\proclaim{Proposition 2.2} Let $p\geq 2$ be an integer. Then there
is a positive constant $c_p$ such that
$$
|1-z|\leq \exp\left(-\R\sum_{k=1}^{p-1} \frac{z^k}{k} +
c_p|z|^p\right) \qquad (z\in\Bbb C).
$$
\endproclaim

\proclaim{Corollary} Suppose that $f$ is a polynomial and $f(0)=1$.
Then
$$
|f(z)|\leq e^{\tilde s_1(f) |z|} \qquad (z\in\Bbb C),
$$
and for $p=2, 3, \dots$ we have
$$
|f(z)|\leq \exp\left(-\R\sum_{k=1}^{p-1} \frac{1}{k}s_k(f) z^k +
c_p\tilde s_p(f)|z|^p\right) \qquad (z\in\Bbb C).
$$
\endproclaim

In \cite{4}, Lindwart and P\'olya proved the following theorem under
the assumption that the sequence converges uniformly on a disk
$D_r$. Fortunately, their proof works in the case of weak
convergence as well.

\proclaim{Theorem 2.3} Suppose that $p$ is a positive integer,
$M>0$, $\langle f_n\rangle$ is a sequence of polynomials, $f_n(0)\ne
0$ for all $n$, $\lim_{n\to\infty}f_n(0)\ne 0$, $f_n\to f$ weakly,
and that $\tilde s_p(f_n)\leq M$ for all $n$. Then $f_n\to f$
strongly, and $f$ is of the form
$$
f(z)=e^{\alpha z^p} g(z)\qquad (z\in\Bbb C),
$$
where $\alpha$  is a constant and $g$ is an entire function of genus
at most $p-1$.
\endproclaim

\demo{Proof} By replacing $\langle f_n\rangle$ with  $\langle
f_n(0)^{-1} f_n\rangle$, we may assume that $f_n(0)=1$ for all $n$.
Since $\langle f_n\rangle$ converges weakly and $\tilde s_p(f_n)\leq
M$ for all $n$, there are positive constants $a$ and $b$ such that
$$
|f_n(z)|\leq a e^{b|z|^p}\qquad (n=1, 2, \dotsc ;\  z\in\Bbb C),
\tag 2.1
$$
by Proposition 2.1 and the corollary to Proposition 2.2. In
particular, the polynomials $f_1, f_2, \dots$ are uniformly bounded
on compact sets in the complex plane, hence Montel's theorem implies
that  $\langle f_n\rangle$ has a subsequence  which converges
strongly to an entire function; and since   $f_n\to f$ weakly, the
Maclaurin series of the entire function and the formal power series
$f$ coincide. Thus $f$ represents an entire function. Furthermore,
the same argument shows that {\it every} subsequence of  $\langle
f_n\rangle$ has a subsequence which converges to $f$ strongly.
Therefore $f_n\to f$ strongly.

Finally, (2.1) implies that
$$
|f(z)|\leq a e^{b|z|^p}\qquad (  z\in\Bbb C),
$$
and we have
$$
\tilde s_p(f)=\lim_{R\to\infty} \tilde s_p(f;
D_R)=\lim_{R\to\infty}\lim_{n\to\infty} \tilde s_p(f_n; D_R)\leq M;
$$
hence the last assertion follows from Hadamard's factorization
theorem.  \qed
\enddemo

For $0\leq c<\infty$ we denote the sector $\left\{re^{i\theta}:
r\geq 0,\  |\theta|\leq\arctan c\right\}$ by $S_c$, which may be
expressed as
$$
S_c=\left\{z\in\Bbb C: |z|\leq \sqrt{1+c^2}\,\R z\right\}.
$$
We denote the closed right half plane $\{z\in\Bbb C : \R z\geq 0\}$
by  $S_{\infty}$. It is clear that $S_c$ is closed under addition
for $0\leq c\leq\infty$. For $p=1, 2, \dots$ and for $0\leq
c\leq\infty$ we put
$$
S_c^{1/p}=\{z\in\Bbb C : z^p\in S_c\}.
$$
Each $S_c^{1/p}$ is a closed subset of the complex plane, and we
have
$$
z\in S_c^{1/p}\setminus\{0\} \Leftrightarrow z^{-1}\in
S_c^{1/p}\setminus\{0\}.
$$

The following theorem plays a crucial role in our proofs of Theorems
1.1 and 1.2.

\proclaim{Theorem 2.4}  Suppose that $p$ is a positive integer,
$\langle f_n\rangle$ is a sequence of polynomials, $f_n \to f$
weakly, $\lim_{n\to\infty}f_n(0)\ne 0$,
  and that $\Cal Z(f_n)\subset
S_{\infty}^{1/p}$ for all $n$. Then $f_n \to f$
  strongly,  $\Cal Z(f)\subset S_{\infty}^{1/p}$, and  $f$ is of the form
  $$
  f(z)=e^{\alpha z^{2p}}g(z) \qquad(z\in\Bbb C), \tag 2.2
  $$
  where $\alpha$ is a constant and $g$ is an entire function of
  genus at most $2p-1$.
\endproclaim

\demo{Proof} First of all, we may assume that $f_n(0)\ne 0$ for all
$n$. To show that $\langle \tilde s_{2p}(f_n)\rangle$ is bounded,
let $c$ be a real number $>1$ and $c_1=2c/(c^2-1)$. We have
$0<c_1<\infty$ and
$$
z\in  S_{\infty}^{1/p} \setminus  S_{c}^{1/p}\  \Rightarrow \
-z^{2p}\in S_{c_1}.
$$
Let $\langle g_n\rangle$ and $\langle h_n\rangle$ be sequences of
polynomials such that $f_n=g_n h_n$, $\Cal Z(g_n)\subset
S_{c}^{1/p}$ and  $\Cal Z(h_n)\subset S_{\infty}^{1/p} \setminus
S_{c}^{1/p}$ for all $n$. Then we have
$$
\tilde s_p(g_n)\leq \sqrt{1+c^2}\, \R s_p(g_n) \leq \sqrt{1+c^2}\,
\R s_p(f_n) \tag 2.3
$$
and
$$
\tilde s_{2p}(h_n)\leq - \sqrt{1+c_1^2}\, \R s_{2p}(h_n) \leq
\sqrt{1+c_1^2}\left| s_{2p}(h_n)\right| \tag 2.4
$$
for all $n$. Since $\langle f_n\rangle$ converges weakly, (2.3) and
Proposition 2.1 imply that $\langle \tilde s_p(g_n) \rangle$ is
bounded, and we have $\tilde s_{2p}(g_n)\leq \left(\tilde
s_p(g_n)\right)^2$ for all $n$. Thus  $\langle \tilde s_{2p}(g_n)
\rangle$ is bounded. Since $s_{2p}(f_n)=s_{2p}(g_n)+s_{2p}(h_n)$ and
$|s_{2p}(g_n)|\leq \tilde s_{2p}(g_n)$ for all $n$, it follows that
$\langle s_{2p}(h_n)\rangle$ is bounded, hence (2.4) implies that
$\langle \tilde s_{2p}(h_n)\rangle$ is bounded. Therefore $\langle
\tilde s_{2p}(f_n)\rangle$ is bounded, and Theorem 2.3 implies that
$f_n\to f$ strongly and $f$ is of the form (2.2). Finally, it is
clear that  $\Cal Z(f)\subset S_{\infty}^{1/p}$, because
$S_{\infty}^{1/p}$ is closed and $\Cal Z(f_n)\subset
S_{\infty}^{1/p}$ for all $n$. \qed
\enddemo

\remark{Remark} This theorem is also proved in \cite{4} under the
assumption that the sequence converges uniformly on some $D_r$.
Unlike in the case of Theorem 2.3, our proof is different from
theirs.
\endremark

A large part of the following theorem is known, but we provide a
detailed proof for completeness as well as for the reader's
convenience.

\proclaim{Theorem 2.5}   Suppose that $\langle f_n \rangle$ is a
weakly convergent  sequence of polynomials,
$\lim_{n\to\infty}f_n(0)\ne 0$,
  and that $\Cal Z(f_n)\subset
S_{\infty}$ for all $n$. Then $\langle f_n \rangle$
  converges strongly, and the limit $f$ is of the form
$$
f(z)=f(0) e^{\alpha z^2 - \beta
z}\prod_j\left(1-\frac{z}{a_j}\right)e^{z/a_j} \qquad (z\in\Bbb C),
$$
where $\alpha, \beta$ are constants, $a_j\in
S_{\infty}\setminus\{0\}$ for all $j$ and $\sum |a_j|^{-2}<\infty$.
Furthermore,
$$
\sum_j \R\! \left(a_j^{-1}\right)\leq \R \beta \tag 2.5
$$
 and $\alpha\geq 0$.
\endproclaim

\demo{Proof} We need only to prove the last two inequalities,
because the remaining assertions are immediate consequences of the
case $p=1$ of Theorem 2.4.

First of all, we may assume that $f_n(0)\ne 0$ for all $n$. Then
Proposition 2.1 implies that
$$
\beta=-\left(\frac{f'}{f}\right)(0)=\lim_{n\to\infty} s_1(f_n)
$$
and
$$
2\alpha - s_2(f)=\left(\frac{f'}{f}\right)'(0)=-\lim_{n\to\infty}
s_2(f_n).
$$

Since $\R\left(a^{-1}\right)\geq 0$ for all $a\in
S_{\infty}\setminus\{0\}$, we have
$$
\sum_j \R\! \left(a_j^{-1}\right)=\sup_{R>0} \R s_1(f; D_R)
$$
and
$$
\align
 \R s_1(f; D_R)&=\R \lim_{n\to\infty}s_1(f_n; D_R)\\
 &=\lim_{n\to\infty}\R s_1(f_n;
 D_R)\\
 &\leq \liminf_{n\to\infty}\,\R s_1(f_n)\\
 &=\R\beta\\
 \endalign
 $$
for every $R>0$. Hence (2.5) holds.

To prove $\alpha\geq 0$, let $c$ be a real number $>1$ and
$c_1=2c/(c^2-1)$. Then we have
$$\align
\left|s_2(f_n; S_c\setminus D_R)\right|&\leq \tilde s_2(f_n;
S_c\setminus D_R)\\
&\leq R^{-1}\tilde s_1(f_n;
S_c)\\
&\leq  R^{-1} \sqrt{1+c^2}\,\R s_1(f_n)\qquad (n=1, 2, \dotsc ; \ R>0),\\
\endalign
$$
and it follows that
$$\align
2\alpha&=s_2(f)-\lim_{n\to\infty} s_2(f_n)\\
&=\lim_{R\to\infty}s_2(f; D_R)-\lim_{n\to\infty} s_2(f_n)\\
&=\lim_{R\to\infty}\lim_{n\to\infty}s_2(f_n; D_R)-\lim_{n\to\infty} s_2(f_n)\\
&=-\lim_{R\to\infty}\lim_{n\to\infty}s_2(f_n;S_{\infty}\setminus D_R)\\
&=-\lim_{R\to\infty}\lim_{n\to\infty}s_2(f_n;S_{c}\setminus D_R)
-\lim_{R\to\infty}\lim_{n\to\infty}s_2(f_n;(S_{\infty}\setminus S_c)\setminus D_R)\\
&=-\lim_{R\to\infty}\lim_{n\to\infty}s_2(f_n;(S_{\infty}\setminus S_c)\setminus D_R).\\
\endalign
$$
If $a\in S_{\infty}\setminus S_c$ and $a\ne 0$, then $-a^{-2}\in
S_{c_1}$, hence $2\alpha\in S_{c_1}$. Since $c>1$ was arbitrary and
$c_1\to 0$ as $c\to\infty$, we conclude that $\alpha\in S_0=[0,
\infty)$. \qed
\enddemo

\remark{Remark} If there is a $c$ such that $0\leq c<\infty$ and
$\Cal Z(f_n)\subset S_c$ for all $n$, then (2.5) and the argument
given in the last paragraph of the   proof show that $\tilde
s_1(f)<\infty$ and $\alpha=0$; hence we have
$$
f(z)=f(0)e^{-\gamma z}\prod_j\left(1-\frac{z}{a_j}\right)\qquad
(z\in\Bbb C),
$$
where $ \gamma=\beta-s_1(f)$. In this case, we have
$$
\gamma=\lim_{R\to\infty}\lim_{n\to\infty}s_1(f_n; S_c\setminus
D_R)\in S_c.
$$
\endremark

\proclaim{Corollary}  Suppose that $\langle f_n \rangle$ is a weakly
convergent  sequence of polynomials, $\lim_{n\to\infty}f_n(0)\ne 0$,
  and that $\Cal Z(f_n)\subset
\Bbb H$ for all $n$. Then $\langle f_n \rangle$
  converges strongly, and the limit is in the P\'olya-Obrechkoff
  class. If the polynomials $f_1, f_2, \dots$ are real, then the
  limit is in the Laguerre-P\'olya class.
\endproclaim

\demo{Proof} If $f$ and $g$ are entire functions and $f(z) = g(iz)$,
then $\Cal Z(f) \subset S_{\infty}$ if and only if $\Cal Z(g)
\subset \Bbb H$, and  $f$ is as in Theorem 2.5 if and only if $g
\in\Cal{PO}$. This proves the first assertion. The second assertion
is obvious. \qed
\enddemo

Now, we are ready to prove Theorem 1.1.

\demo{Proof of Theorem 1.1} First of all, it is trivial to see that
if $\langle f_n\rangle$ converges weakly and $\langle
f_n^{(k)}\rangle$ converges strongly for some $k$, then the original
sequence converges strongly.

 Let $k$ be such that $ \lim_{n\to\infty}f_n^{(k)}(0)\ne 0$. Since $f_1, f_2,
\dots$ are polynomials, and since $\Cal Z(f_n)\subset \Bbb H$ for
all $n$, it follows from the Gauss-Lucas theorem that $\Cal
Z(f_n^{(k)})\subset \Bbb H$ for all $n$. We have
$\lim_{n\to\infty}f_n^{(k)}(0)\ne 0$ and it is obvious that $\langle
f_n^{(k)}\rangle$ converges weakly. Hence the corollary to Theorem
2.5 implies that $\langle f_n^{(k)}\rangle$ converges strongly, and
we conclude that $f_n\to f$ strongly.

It remains to show that $f\in\Cal{PO}$. If $t\in\Bbb R$ and
$g(z)=f(z+t)$, then it is easy to see that $f\in\Cal {PO}$ if and
only if $g\in\Cal {PO}$. Since $f^{(k)}(0)\ne 0$, there is a
$t\in\Bbb R$ such that $f(t)\ne 0$. If we put $g(z)=f(z+t)$ and
$g_n(z)=f_n(z+t)$, then $g_n\to g$ strongly, $g(0)\ne 0$ and $\Cal
Z(g_n)\subset\Bbb H$ for all $n$, hence the corollary to Theorem 2.5
implies that $g\in\Cal{PO}$, as desired. \qed
\enddemo

In order to prove Theorem 1.2, we need a technical result.

\proclaim{Proposition 2.6} Let $N$ be a positive integer. Then there
is a finite set $\Cal Q$ of positive even integers having the
following property: If $z_1, \dots, z_N\in\Bbb C$, then there is
some $q\in \Cal Q$ such that $z_1^{q}, \dots, z_N^{q}\in
S_{\infty}$.
\endproclaim

\demo{Proof} The unit circle $\Bbb T = \{\zeta\in\Bbb C :
|\zeta|=1\}$ is a compact set in $\Bbb C$ and is a multiplicative
group. For $q=2,4,   \dots$ we put
$$
U_q=\{(\zeta_1, \dots, \zeta_N)\in\Bbb T^N: \zeta_1^{q}, \dots,
\zeta_N^{q}\in\, \text{int}\, S_{\infty}\}.
$$
It is clear that $U_2, U_4, \dots$ are open subsets of the compact
metric space $\Bbb T^N$. If $(\zeta_1, \dots, \zeta_N)\in\Bbb T^N$,
then
 $\zeta_1^{q}, \dots,
\zeta_N^{q}\in \text{int}\, S_{\infty}$ for some positive even
integer $q$. In fact, if $g$ is an element of a sequentially compact
topological group, any neighborhood of the unit element contains
$g^n$ for infinitely many positive integers $n$.
  We have shown that $\{ U_{q} : q=2, 4, \dots\}$ is an
open cover of $\Bbb T^N$. Hence there is a finite set $\Cal Q$ of
positive even integers having the property. \qed
\enddemo

 If $t\in\Bbb R$  and  $g(z)=f(z+t)$, then (as in the proof of Theorem 1.1) $f\in\Cal{LP}^*$ if and only
if $g\in\Cal{LP}^*$.

\demo{Proof of  Theorem 1.2}   We first consider the case where
$\lim_{n\to\infty}f_n(0)\ne 0$. In this case,  we may assume that
$f_n(0)=1$ for all $n$.  Let $N$ be a positive integer such that
$N(f_n;\Bbb C\setminus\Bbb R)\leq N$ for all $n$. Then there are
sequences $\langle P_n\rangle$ and $\langle g_n\rangle$ of real
polynomials such that $f_n=P_n g_n$, $P_n(0)=g_n(0)=1$, $\Cal
Z(P_n)\cap \Bbb R=\emptyset$, $\deg P_n\leq N$ and $\Cal
Z(g_n)\subset \Bbb R$ for all $n$. Let $\Cal Q$ be as in Proposition
2.6. Since each $P_n$ has at most $N$ zeros, it follows that for
every $n$ there is a $q\in \Cal Q$ such that $\Cal Z (P_n)\subset
S_{\infty}^{1/q}$. Since every $q\in \Cal Q$ is even and $\Cal
Z(g_n)\subset \Bbb R$ for all $n$, we have $\Cal Z(g_n)\subset
S_{\infty}^{1/q}$ for all $n$ and for all $q\in \Cal Q$.
 Consequently, for every $n$ there is a $q\in \Cal Q$
such that $\Cal Z (f_n)\subset S_{\infty}^{1/q}$; and it follows
from Theorem 2.4 that $\langle f_n \rangle$ is the union of a finite
number ($\leq |\Cal Q|$) of strongly convergent subsequences.
Therefore $f_n\to f$ strongly.

Since $f_n\to f$ strongly and $f_n(0)=1$ for all $n$, there is a
positive constant $r$ such that $\Cal Z(f_n)\cap D_r=\emptyset$ for
all $n$, and it follows that $\Cal Z(P_n)\cap D_r=\emptyset$ for all
$n$. Since $\deg P_n\leq N$ and $P_n(0)=1$ for all $n$, the
coefficients of the polynomials $P_1, P_2, \dots$ are uniformly
bounded. Hence there is a strictly increasing sequence $\langle
n(l)\rangle$ of positive integers such that $\langle
P_{n(l)}\rangle$ converges strongly to a real polynomial $P$, and it
follows that $\langle g_{n(l)}\rangle$ converges uniformly on a disk
centered at the origin. Since $g_n(0)=1$ and $\Cal Z(g_n)\subset\Bbb
R$ ($\subset\Bbb H$) for all $n$, and since $g_1, g_2, \dots$ are
real polynomials, the corollary to Theorem 2.5 implies that $\langle
g_{n(l)}\rangle$ converges strongly to an entire function
$g\in\Cal{LP}$. Therefore $f=Pg\in\Cal{LP}^*$.

In the general case, there is an integer $k$ such that
$\lim_{n\to\infty}f_n^{(k)}(0)\ne 0$. Since $f_1, f_2, \dots$ are
real  polynomials, Rolle's theorem implies that $N(f_n^{(k)};\Bbb
C\setminus\Bbb R)\leq N(f_n;\Bbb C\setminus\Bbb R)$ for all $n$.
Hence $\langle f_n^{(k)}\rangle$ converges strongly, and it follows
that $\langle f_n \rangle$ converges strongly. Since $f^{(k)}(0)\ne
0$, there is some $t\in\Bbb R$ such that $f(t)\ne 0$, and the same
argument as in the proof of Theorem 1.1  shows that $f
\in\Cal{LP}^*$. \qed
\enddemo

\heading 3. Zeros of Appell and Jensen Polynomials
\endheading

Let $f$ be a formal power series given by
$$
f(z)=\sum_{k=0}^{\infty} a_k z^k.
$$
 For $n=0, 1, 2, \dots$ the $n$-th
{\it Appell}  and {\it Jensen} polynomials $A{(f,n)}$ and $J{(f,n)}$
of $f$ are defined by
$$
A{(f,n)}(z)=\sum_{k=0}^n\frac{n!}{(n-k)!} a_k z^{n-k}\ \ \text{and}\
\  J{(f,n)}(z)=\sum_{k=0}^n\frac{n!}{(n-k)!} a_k z^{k},
$$
respectively.  We have
$$
A{(f,n)}(z)=z^n J{(f,n)}\left(\frac{1}{z}\right),
$$
and it follows that $N(A{(f,n)}; \Bbb C\setminus \Bbb R)=N(J{(f,n)};
\Bbb C\setminus \Bbb R)$ and $N(A{(f,n)}; \Bbb C\setminus
S_0)=N(J{(f,n)}; \Bbb C\setminus S_0)$ for all $n$.
 As in \cite{1}, we define the sequence   $\langle f^*_n\rangle$ of polynomials by
$$
f^*_n(z)=J{(f,n)}\left(\frac{z}{n}\right) \qquad (n=1,2, \dots),
$$
which may be called the {\it Jensen sequence} of $f$. It is easy to
see that $f^*_n\to f$ weakly, in the general case; that $f^*_n\to f$
uniformly on a disk $D_r$ if and only if the radius of convergence
is positive; and that $f^*_n\to f$ strongly if and only if $f$
represents an entire function. It is also easy to see that if
$f_k\to f$ weakly as $k\to\infty$, then $A(f_k, n)\to A(f, n)$ and
$J(f_k, n)\to J(f, n)$ strongly as $k\to\infty$ for every $n$. The
following two propositions are known, and will be proved at the end
of this section.

\proclaim{Proposition 3.1} If $f$ is a real polynomial, then
$$
N(J(f,n);\Bbb C\setminus\Bbb R)\leq N(f; \Bbb C\setminus\Bbb R)
$$
and
$$
 N(J(f,n);\Bbb C\setminus S_0)\leq N(f; \Bbb C\setminus S_0)
$$
for all $n$.
\endproclaim

\remark{Remark} Since  $N(A{(f,n)}; \Bbb C\setminus \Bbb
R)=N(J{(f,n)}; \Bbb C\setminus \Bbb R)$ and $N(A{(f,n)}; \Bbb
C\setminus S_0)=N(J{(f,n)}; \Bbb C\setminus S_0)$ for all $n$, the
term $J(f,n)$ may be replaced with $A(f,n)$.
\endremark

\proclaim{Proposition 3.2} If $f$ is a polynomial and $\Cal
Z(f)\subset \Bbb H$, then   $\Cal Z(J(f,n))\subset \Bbb H$ for all
$n$.
\endproclaim

In 1914, P\'olya and Schur proved that a formal power series $f$
with real coefficients represents an entire function in $\Cal{LP}$
if and only if  $\Cal Z(J(f,n))\subset \Bbb R$ for all $n$ \cite{7,
\S 6 Satz IV}, and later Benz generalized the result \cite{1, Satz
2}. Theorem 1.1 immediately implies Benz's theorem.

\proclaim{Theorem 3.3} Suppose that $f$ is a formal power series.
Then the following are equivalent:
 \roster
 \item $f$ represents an entire function in $\Cal{PO}$.
\item There  is a
sequence $\langle f_k\rangle$ of polynomials such that $f_k\to f$
weakly and $\Cal Z(f_k)\subset\Bbb H$ for all $k$.
\item $\Cal Z(J(f,n))\subset\Bbb H$ for all $n$.
\endroster

\endproclaim
\demo{Proof}
 If  $f$ represents an entire function in $\Cal{PO}$,
then there is a sequence $\langle f_k\rangle$ of polynomials such
that $f_k\to f$ strongly (hence weakly) and $\Cal Z(f_k)\subset\Bbb
H$ for all $k$.

Next, suppose that  $\langle f_k\rangle$ is a sequence of
polynomials such that $f_k\to f$ weakly and $\Cal Z(f_k)\subset\Bbb
H$ for all $k$.  Let $n$ be arbitrary. Then $J(f_k, n)\to J(f, n)$
strongly as $k\to\infty$, and Proposition 3.2 implies that $\Cal
Z(J(f_k,n))\subset \Bbb H$ for all $k$. Therefore we have $\Cal
Z(J(f,n))\subset\Bbb H$.

Finally, suppose that $\Cal Z(J(f,n))\subset\Bbb H$ for all $n$. If
every coefficient of $f$ is $0$, then it is obvious that
$f\in\Cal{PO}$, otherwise the Jensen sequence $\langle f^*_n\rangle$
satisfies the conditions of Theorem 1.1, and we conclude that $f$
represents an entire function in $\Cal{PO}$.
 \qed
\enddemo

\remark{Remarks} (i) If (2) holds, then Theorem 1.1 implies that the
sequence converges strongly, unless $f$ is identically equal to $0$.
(ii) Benz proved the implication (3) $\Rightarrow$ (1) by showing
that the $f_n^*\to f$ strongly, but the argument does not imply
Theorem 1.1.
\endremark

 The Laguerre-P\'olya class has a subclass which
plays an important role in the theory of multiplier sequences
\cite{7}. Let $\Cal{LP}( S_0)$ be the class of real entire functions
$f$ which are of the form
$$
f(z)= c z^m e^{-\alpha z}\prod_j\left(1-\frac{z}{a_j}\right)\qquad
(z\in\Bbb C),
$$
where $c\in\Bbb R$,   $m$ is a non-negative integer, $\alpha\geq 0$,
$a_1, a_2, \dots >0$ and $\sum_j a_j^{-1}<\infty$. Likewise,
$\Cal{LP}( S_0)^*$ will denote the class of real entire functions of
the form $Pg$ where $P$ is a real polynomial and $g\in \Cal{LP}(
S_0)$. It is clear that
$$
\Cal{LP}( S_0)=\left\{f\in\Cal{LP}( S_0)^* : \Cal Z(f)\subset [0,
\infty)\right\},
$$
and it is known that $f\in\Cal{LP}( S_0)^*$ if and only if there is
a sequence $\langle f_n\rangle$ of real polynomials such that
$f_n\to f$ strongly and $\langle N(f_n; \Bbb C\setminus S_0)\rangle$
is bounded \cite{5, \S 4}. Since $N(f; \Bbb C\setminus \Bbb R)\leq
N(f; \Bbb C\setminus S_0)$ for every entire function $f$, Theorem
1.2 implies the following:

\proclaim{Theorem 3.4}  Suppose that  $\langle f_n\rangle$ is a
sequence of real polynomials, $f_n\to f$ weakly, some coefficient of
$f$ is not $0$, and that $\langle N(f_n; \Bbb C\setminus
S_0)\rangle$ is bounded. Then $f$ represents an entire function in
$\Cal{LP}( S_0)^*$ and $f_n\to f$ strongly.

\endproclaim

\remark{Remark} This theorem may be proved directly by the same
argument (with $S_0$ instead of $\Bbb R$) as in the proof of Theorem
1.2. See the remark after the proof of Theorem 2.5.
\endremark

In the following two theorems, P\'olya proved (1) and (2),  proved
(3) under the additional  assumption that the radius of convergence
of $f$ is positive  \cite{5, \S 4; 6, Satz II}, and conjectured that
the
  assumption could be unnecessary \cite{5, p.246 footnote
$^*)$}. It seems that no progress has
  been made about the conjecture.

  \proclaim{Theorem 3.5} Let $f$ be a formal power series with real
coefficients, and for  $n=0, 1, 2, \dots$ let $N_n=N(A(f,n); \Bbb
C\setminus S_0)$. Then
 \roster
 \item $0=N_0\leq N_1\leq N_2 \leq\dotsb$;
\item if $f$ represents an entire function in
$\Cal{LP}( S_0)^*$ and $N=N(f; \Bbb C\setminus S_0)$, then $N_n\leq
N$ for all $n$, and  $N_n= N$ for all sufficiently large $n$; and
\item if  $\langle N_n\rangle$
 is bounded, then  $f$ represents an entire function in
$\Cal{LP}(S_0)^*$.
\endroster
\endproclaim

\proclaim{Theorem 3.6} Let $f$ be a formal power series with real
coefficients, and for $n=0, 1, 2, \dots$ let $N_n=N(A{(f,n)};\Bbb
C\setminus\Bbb R)$. Then
 \roster
 \item $0=N_0=N_1\leq N_2\leq N_3\leq\dotsb$;
\item if $f$ represents an entire function in
$\Cal{LP}^*$ and $N= N(f;\Bbb  C\setminus\Bbb R)$, then $N_n\leq  N$
for all $n$, and  $N_n= N$ for all sufficiently large $n$; and
\item if  $\langle N_n\rangle$
 is bounded, then  $f$ represents an entire function in
$\Cal{LP}^*$.
\endroster
\endproclaim

Since the proofs of these theorems are almost identical, we prove
the first one only.

\demo{Proof of Theorem 3.5} The equality in (1) is obvious, and the
remaining inequalities follow from Rolle's theorem, because we have
$A(f,n)'=nA(f,n-1)$  for every positive integer $n$.

To prove (2), suppose that $f\in\Cal {LP}( S_0)^*$ and  $N=N(f;\Bbb
C\setminus S_0)$. Then there is a sequence $\langle f_k\rangle$ of
real polynomials such that $f_k\to f$ strongly as $k\to\infty$ and
$N(f_k; \Bbb C\setminus S_0)=N$ for all $k$. Let $n$ be an arbitrary
positive integer. Then $J(f_k, n)\to J(f, n)$ strongly as
$k\to\infty$, hence
$$
N_n=N(J{(f,n)};\Bbb C\setminus S_0)\leq N(J{(f_k,n)}; \Bbb
C\setminus S_0)
$$
for all sufficiently large $k$. On the other hand, Proposition 3.1
implies that
$$
N(J(f_k,n); \Bbb C\setminus S_0)\leq N( f_k ;\Bbb C\setminus S_0)=N
$$
for all $k$. Hence $N_n\leq N$.

 Since $f$ represents an entire function, it follows that the Jensen
 sequence $\langle f^*_n\rangle$ converges to
$ f$ strongly, and we have $N(f^*_n;\Bbb C\setminus S_0)=N_n$ for
all $n$. Hence we have
$$
N\leq N_n
$$
for all sufficiently large $n$, and (2) is proved.

Finally, (3) is proved by the same argument as in the proof of the
implication (3) $\Rightarrow$ (1) in  Theorem 3.3, except that
Theorem 3.4 is applied. \qed
\enddemo

\demo{Proof of Proposition 3.1} Suppose that $f$ is a real
polynomial and $N=N(f; \Bbb C\setminus\Bbb R)$. If $f$ has no real
zeros, so that $N=\deg f$, then
$$
N(J{(f,n)}; \Bbb C\setminus\Bbb R )\leq\deg J{(f,n)}\leq\deg f =N
$$
for every non-negative integer $n$.  On the other hand, it is
trivial to see that if $a$ is a constant and $f(z)=(z-a)g(z)$, then
$$
A(f,n)=A(g,n)' - aA(g,n)
$$
for all $n$. Since $N(J{(f,n)}; \Bbb C\setminus\Bbb R )=N(A{(f,n)};
\Bbb C\setminus\Bbb R )$ for all $n$, the first inequality follows
from an induction (on the number of real zeros of $f$) based on the
Hermite-Poulain theorem, which states that if $h$ is a real
polynomial and $b$ is a real constant, then
$$
N(h'-b h; \Bbb C\setminus\Bbb R)\leq N(h; \Bbb C\setminus\Bbb R).
$$

The second inequality is proved in the same way except that the
induction is on the number of non-negative real zeros of $f$ and the
corresponding Hermite-Poulain theorem states that if $h$ is a real
polynomial and $b\geq 0$, then
$$
N(h'-b h; \Bbb C\setminus S_0)\leq N(h; \Bbb C\setminus S_0). \qed
$$
\enddemo

\demo{Proof of Proposition 3.2} Suppose that $a\in \Bbb H$,
$f(z)=(z-a)g(z)$, $n$ is a non-negative integer and $\Cal
Z(J(g,n))\subset\Bbb H$. Then $A(g,n)$ is a polynomial of degree
$\leq n$,  $\Cal Z(A(g,n))\subset -\Bbb H$ and
$$
A(f,n)=A(g,n)' - aA(g,n).
$$
Let $d=\deg A(g,n)$ and let $a_1, \dots, a_d$ be the zeros of $A(g,
n)$. Then $\I a_j\leq 0$ for all $j$, and we have
$$
\frac{A(g,n)'(z)}{A(g,n)(z)}=\sum_{j=1}^d\frac{1}{z-a_j}\qquad (z\ne
a_1, \dots, a_d).
$$
If $\I z>0$, then $A(g,n)(z)\ne 0$, and we have
$$
\I \frac{A(f,n)(z)}{A(g,n)(z)}= -\sum_{j=1}^d\frac{\I
(z-a_j)}{|z-a_j|^2} - \I a <0.
$$
Thus  $\Cal Z(A(f,n))\subset -\Bbb H$, and we conclude that   $\Cal
Z(J(f,n))\subset \Bbb H$.

Now it is clear that the result follows from an induction on the
degree of $f$. \qed

\enddemo

\enddocument